\documentclass{article}

\usepackage{fullpage}
\usepackage{amsmath}
\usepackage{amssymb}
\usepackage[all]{xy}
\numberwithin{equation}{section}
\newcommand{\bA}{\mathbb{A}}

\newcommand{\bC}{\mathbb{C}}

\newcommand{\bP}{\mathbb{P}}
\newcommand{\bQ}{\mathbb{Q}}

\newcommand{\cF}{\mathcal{F}}

\newcommand{\cO}{\mathcal{O}}
\newcommand{\cT}{\mathcal{T}}
\newcommand{\cV}{\mathcal{V}}
\newcommand{\cU}{\mathcal{U}}

\newcommand{\wt}{\widetilde}
\newcommand{\ol}{\overline}
\newcommand{\nd}{\noindent}

\newcommand{\Spec}{\text{Spec}\,}

\newcommand{\D}{\Delta}
\newcommand{\proof}{\textbf{Proof. }}

\makeatletter
\renewcommand{\section}{\@startsection{section}{1}{0mm}{1.5\baselineskip}{0.5\baselineskip}{\large\bf\center}}
\renewcommand{\subsection}{\@startsection{subsection}{2}{5mm}{0.5\baselineskip}{-0.5em}{\bf }}
\renewcommand{\subsubsection}{\@startsection{subsubsection}{3}{5mm}{0.5\baselineskip}{-0.5em}{\bf}}
\renewcommand\thesection{\@arabic\c@section}
\renewcommand\thesubsection{\thesection.\@arabic\c@subsection}
\renewcommand\thesubsubsection{\thesubsection.\@arabic\c@subsubsection}
\renewcommand\theequation{\thesection.\@arabic\c@equation}
\renewcommand\appendix{\par\setcounter{section}{0}\setcounter{subsection}{0}\gdef\thesection{\@Alph\c@section}}
\makeatother

\begin{document}

\begin{center}
{\Large \bf Some remarks on Landau-Ginzburg potentials\\ for odd-dimensional quadrics}

\vspace{10pt}

Vassily Gorbounov, Maxim Smirnov

\vspace{10pt}

\end{center}

\begin{abstract}
\nd We study the possibility of constructing a Frobenius manifold for the standard Landau-Ginzburg model of odd-dimensional quadrics $Q_{2n+1}$ and matching it with the Frobenius manifold attached to the quantum cohomology of these quadrics. Namely, we show that the initial conditions of the quantum cohomology Frobenius manifold of the quadric can be obtained from the suitably modified standard Landau-Ginzburg model. 

\end{abstract}


\section{Introduction}

The idea of the Landau-Ginzburg model, LG model from now on, of a variety originated in the work of physicists. In studying conformal field theories over projective hypersurfaces, the zeros of a certain function, they noticed that the function itself can be used to calculate various quantities in the theory of the hypersurface. This revealed a mysterious relation between smooth projective varieties and functions with isolated singularities, the superpotentials, as they are called in the literature. We cite, out of many references, the work of D.~Gepner,  E.~Witten, and T.~Eguchi and his school as being very influential for mathematicians \cite{Ge}, \cite{Wi}, \cite{EgHoXi}.

The LG models appeared in the mathematics literature for the first time probably in the work of A.~Givental \cite{Gi}. He has introduced the notion of the quantum differential equation associated to a smooth projective variety, and expressed its solutions for some classes of Fano varieties as oscillating integrals of a certain function, the "LG model" of these Fano varieties.

Later the concept of the LG model became a part of the Homological Mirror Symmetry Conjecture, which is due to M.~Kontsevich.

Despite a large volume of serious work of many people there is no commonly accepted definition of the LG model for a given Fano variety. There is a number of examples of a LG model for Fano varieties but they are examples each in its own sense.

In developing the theory of Frobenius manifolds, B.~Dubrovin used two main sources of examples. The first is the quantum cohomology of a smooth projective variety, just a variety from now on, and the other is the work of K.~Saito \cite{Sai}, an important development of singularity theory, which produces a non-trivial example of a Frobenius manifold associated to a polynomial superpotential with an isolated singularity. This suggests another natural approach to finding the LG model for a variety $X$: one could look for a function on a non-compact manifold such that some generalization of the K.~Saito construction would produce a Frobenius manifold isomorphic to the one the quantum cohomology of $X$ produces. This approach to the LG model was detailed in the book of Yu.~Manin \cite{Ma}. It was actually made to work for the weighted projective spaces in the papers of A.~Douai, C.~Sabbah, and E.~Mann \cite{DoSa1}, \cite{DoSa2}, \cite{M}.

In this note we expand \cite{DoSa2} by considering varieties which are not toric. Namely, we take one of the known in the literature Laurent polynomials which is believed to be the LG model for an odd-dimensional quadric (see \cite{HoVa}, \cite{Pr}), and investigate the possibility of making it an example of a LG model to the quadric in the sense of the Frobenius manifold structure. From now on we will refer to this potential as "the standard potential".

We attempted to formulate a rule for finding such a potential in general. In the case of the quadric as well as in some other examples available in the literature the partial derivatives of the standard potential reproduce partially the relations generated by the multiplications by the generator of the second cohomology. But here a problem arises: the standard potential may have fewer critical points than the rank of the cohomology of the variety it supposedly models, as is observed in the quadric setting. Therefore this LG potential has the wrong Jacobi ring and in terms of the corresponding Frobenius manifolds can't be the correct LG model. This problem for the LG model of the Grassmannians was pointed out in the paper \cite{EgHoXi} a long time ago. In the same paper it was suggested how to remedy the situation. One needs to search for a partial compactification of the domain of the standard potential to gain the required number of the critical points. There is no procedure for doing this in general at least as far as we know, but in the examples we have looked at, including the quadric, it was always possible to find such a compactification by trial and error. The compactification of the domain changes the superpotential and it is no longer a Laurent polynomial.

We follow the generalization due to C.~Sabbah and A.~Douai of the K.~Saito work. This reduces the problem of construction of a Frobenius manifold to studying the Gauss-Manin system associated to the potential. The first step is to solve the Birkhoff problem for the Gauss-Manin system at one point of the parameter space for the  potential, to make sure it matches the initial conditions of the quantum cohomology Frobenius manifold of the appropriate variety. In the case of the quadric this can be done but with bit of an effort. With the kind help of C.~Sabbah and A.~Nemethi we check a certain property for our potential, the so called cohomological tameness, which helped a lot with solving the Birkhoff problem.

For the final step, we need to build the Frobenius manifold in the neighbourhood of that one point in the parameter space which would automatically match the Frobenius manifold of the quadric. This amounts to finding a good deformation of our potential, and this in turn requires one to have a good hold on the behaviour of the potential at the infinity. It is not clear at the present time what the best method is to handle this purely analytic problem in its most general setting. There are, however, some sufficiency conditions -- for example the so called $M$-tameness of the potential -- that guarantee such behaviour. In the case of Laurent polynomial this property can be checked and is used in \cite{DoSa1} \cite{DoSa2}. For a general regular function it appears to be a difficult analytical problem. We fail to check the M-tameness for our potential, however it seems interesting to us that proving the existence of the LG model hinges on the delicate analytical properties of the potential.

\smallskip

Putting it all together, in this note we have proved that the standard potential for the quadric after a suitable modification defines the initial conditions of the quantum cohomology Frobenius manifold of the quadric.

\smallskip

Throughout the whole text we restrict ourselves to the case of a three-dimensional quadric $Q_3$ for the sake of making the exposition compact. In fact, all the results hold for an arbitrary odd-dimensional quadric.

\medskip

\textbf{Remark.} After the first version of this preprint was posted on the arxiv an interesting preprint by C.~Pech and K.~Rietsch \cite{PeRi} appeared where it is shown that the superpotential considered in this paper for an odd-dimensional quadric is a particular case of a general construction proposed by K.~Rietsch \cite{Ri} for a homogeneous space $G/P$. Some results of Section \ref{Sec.: Frobenius for Quadric} are used in \cite{PeRi} to establish a part of a conjecture made in \cite{Ri} (see \textit{loc.cit.} for details). It would be interesting to further bridge the work done in this preprint with the approach initiated by K.~Rietsch.

\medskip

\textbf{Acknowledgements.} First, we would like to thank Yu.~I.~Manin for initiating this project and sharing generously his ideas. We are indebted to C.~Sabbah and A.~Nemethi for the proof of Lemma \ref{SubSubSec.: Lemma Tameness}, and help with general questions about the subject of this paper. Special thanks go to P.~Bressler and to E.~Shinder for various discussions related to this work.

Both authors benefited a lot from the excellent research environment of the MPIM, Bonn. Some part of this work was done while the second author visited the IHES, Paris whose hospitality is gratefully acknowledged.

\section{Background and notation}
\label{Sec.: Background and notations}

\subsection{Quantum cohomology.}
\label{SubSubSec.: Quantum cohomology}

Here we very briefly recall some facts about quantum cohomology. For the general theory of Frobenius manifolds we refer to \cite{Ma}, \cite{He}, and references therein.

Let $X$ be a smooth projective complex algebraic variety and $QH(X)$ its big quantum cohomology. If $\D_0, \dots, \D_r$ is a graded basis and $x_0, \dots, x_r$ dual coordinates, then the quantum product is defined as
\begin{align}\label{Eq.: Quantum Product}
\D_i \circ \D_j = \sum_{k,l} \Phi_{ijk} g^{kl} \D_l,
\end{align}
where $\Phi$ is the Gromov-Witten potential, $\Phi_{ijk}=\frac{\partial^3\Phi}{\partial x_i \partial x_j \partial x_k}$, and $g$ is the Poincar\'e pairing on $H=H^*(X,\mathbb{C})$. The full structure of the quantum cohomology of $X$ endows $H$ with a structure of a Frobenius manifold. In fact, one needs to work with formal manifolds because $\Phi$ is not known to be convergent.

Disregarding convergence issues one can think of $QH(X)$ as a family of multiplications on $H$ parametrized by $H$ itself, i.e. it is a multiplication on $\mathcal{T}_H$, where we consider $H$ as a complex manifold. Poincar\'e pairing on $H$ defines a constant pairing on $\mathcal{T}_H$ which is multiplication invariant and flat.

Under small quantum cohomology one means the restriction of the above picture to $H^2(X, \mathbb{C}) \subset H$. In terms of coordinates it means that we reduce all our formulas modulo an ideal generated by coordinates dual to $\D_i$'s not lying in $H^2(X, \mathbb{C})$.

\subsubsection{Spectral cover.}
\label{SubSubSec.: Spectral Cover for QH}

Assume that $H=H^*(X, \mathbb{C})$ is of dimension one in each even degree and zero otherwise. In this situation one can consider an algebraic torus $\textbf{T} \subset H^t$, which is a locally closed subvariety of $H^{t}$. Here $H^t$ is the dual of $H$.

Namely, let $\D_0, \dots, \D_r$ be a graded basis of $H$, such that $\D_0$ is the identity element, and consider
$$
H^t = \text{Spec}\,(S^{\bullet}(H)) = \text{Spec}\,(\mathbb{C}[\D_0, \dots , \D_r]).
$$
In $H^t$ we have an affine subspace $\{\D_0=1\} \simeq \text{Spec}\,(\mathbb{C}[\D_1, \dots , \D_r])$, and inside this affine subspace we have the torus $\textbf{T}=\text{Spec}\,(\mathbb{C}[\D_1^{\pm 1}, \dots , \D_r^{\pm 1}])$. This torus does not depend on the choice of $\D_1, \dots, \D_r$ and will play an important role in construction of LG models.

Equations that define the spectral cover $\text{Spec}\,(QH(X))$ as a subvariety of  $H \times H^t$ are given just by the multiplication table. One of the equations is always $\D_0=1$.  Hence, the spectral cover always lives inside the affine space $\{\D_0=1\} \simeq \text{Spec}\,(\mathbb{C}[\D_1, \dots , \D_r])$. 

One can summarize it in the diagram
\begin{align*}
\xymatrix{
\text{Spec}\,(QH(X)) \ar[r]^<<<<<i \ar[rd]& H \times H^t  \ar[d]  & \ar[l]_<<<<j\ar[ld] H \times \textbf{T}\\
                                                   & H  
}
\end{align*}
where $i$ and $j$ are embeddings.

In some cases $\text{Spec}\,(QH(X))$ lies in $H \times \textbf{T}$. For example, it is true for projective spaces (at least in the small quantum cohomology). It is not true for the case of odd-dimensional quadrics that we consider.

\subsection{Landau-Ginzburg models.}
\label{SubSec.: LG models}

Let $X$ be a Fano variety and $QH(X)$ its quantum cohomology. 

A Saito's framework (see \cite[III.8]{Ma}) is called a Landau-Ginzburg model for $X$ iff it is isomorphic to $QH(X)$ as a Frobenius manifold. 

Consider a pair $(U, f)$ consisting of a complex smooth affine variety and a regular function on it. It is called a Landau-Ginzburg model for $X$ iff there exists a deformation of $(U, f)$ and a Saito's framework attached to it which is isomorphic to $QH(X)$ as a Frobenius manifold. 

These definitions are quite restrictive. One can relax both of them as follows. We require the existence of a point in $QH(X)$ and a point in Saito's framework such that the germs of Frobenius manifolds at these points are isomorphic.

\subsection{Gauss-Manin systems.}
\label{SubSec.: GM system}

Here we will give a brief account on Gauss-Manin systems, Brieskorn lattices, higher residue pairings and Birkhoff problems. Mainly we will be setting up notation that we will use in Section \ref{Sec.: Frobenius for Quadric} and refer to \cite{Do1}, \cite{Do2} for details.

\smallskip

Consider a projective line with a chosen coordinate $\theta$. We will denote it by $\mathbb{P}^1_{\theta}$. Let $\mathbb{P}^1_{\theta}=U_0\cup U_{\infty}$ be the standard open cover, i.e. $U_0=\mathbb{A}^1_{\theta}$, $U_1=\mathbb{A}^1_{\theta^{-1}}$, and $W=U_0 \cap U_{\infty}=\mathbb{A}^1_{\theta}-\{0\}$. Here $\mathbb{A}^1_{\theta}$ stands for $\text{Spec}(\mathbb{C}[\theta])$ and $\{0\}$ for $\{\theta=0\}$. Sometimes we will write $\tau$ for $\theta^{-1}$. We will use this notation when working with Gauss-Manin systems throughout the article.

\subsubsection{Definition.}

Let $X$ be a smooth affine variety of dimension $n$ with a regular function $h$ on it. The main example relevant for mirror symmetry: $X=(\mathbf{G}_m)^n$ and $h$ is a Laurent polynomial. To such a function $h$ one can attach its \textit{Gauss-Manin system}
\begin{align}\notag
G=\Omega^n(X)[\theta, \theta^{-1}]/(\theta d-dh\wedge) \Omega^{n-1}(X)[\theta, \theta^{-1}],
\end{align}
which is a free $\mathbb{C}[\theta, \theta^{-1}]$-module of finite rank with a flat connection $\nabla$ defined as follows. Let $\sum_i \omega_i\theta^i$ be a representative of some class $\gamma \in G$, i.e. $\gamma=[\,\sum_i \omega_i\theta^i\,]$. Then
\begin{align}\notag
\theta^2\nabla_{\frac{\partial}{\partial \theta}}(\gamma)=\left[\sum_i h\omega_i\theta^i+\sum_i i \omega_i\theta^{i+1}\right],
\end{align}
where the brackets $[\,\,\,\,]$ denote taking class in $G$.

It is a general fact that $(G, \nabla)$ always has a regular singularity at $\theta=\infty$ and possibly irregular singularity at $\theta=0$ of rank 2.

\subsubsection{Brieskorn lattice.}

At $\theta=0$ the Gauss-Manin system $(G, \nabla)$ has a natural lattice 

\begin{align}\label{Eq.: Brieskorn Lattice}
G_0=\Omega^n(X)[\theta]/(\theta d-dh\wedge) \Omega^{n-1}(X)[\theta],
\end{align}
the \textit{Brieskorn lattice} of $h$. This means that $G_0$ is a $\bC[\theta]$-module such that $G_0 \otimes_{\bC[\theta]}\bC[\theta, \theta^{-1}] \simeq G $. 

The connection $\nabla$ naturally restricted (or extended) to $G_0$ has a pole of order 2 at $\theta=0$, i.e.
\begin{align*}
&\theta^2\nabla_{\frac{\partial}{\partial \theta}}(G_0) \subset G_0.
\end{align*}
Abusing notation we will denote this meromorphic connection again by $\nabla$.

In general $G_0$ does not have torsion but may be not finitely generated. If $h$ is \textit{cohomologically tame} (see Section \ref{SubSubSec.: Tameness} and \cite{Sa2}), then it is finitely generated and free. Therefore, in the cohomologically tame case the pair $(G_0, \nabla)$ gives a meromorphic extension of $(G, \nabla)$ to $\bA^1_{\theta}$ with a pole of order 2 at the origin.

\subsubsection{Extension to $\mathbb{P}^1_{\theta}$.} 
\label{SubSubSec.: Extension to P1}

The aim is to find an extension of $(G_0, \nabla)$ to a free $\cO_{\bP^1_{\theta}}$-module with a meromorphic connection on $\bP^1_{\theta}$ with a pole of order less or equal to 1 at infinity. We will denote such an extension $(\cF, \nabla)$. This type of question is known as \textit{the Birkhoff problem} (cf. \cite[Ch. 4]{Sa1}).

More concretely it means that we need to find a $\bC[\theta^{-1}]$-module $G_{\infty} \subset G$ such that: 

\begin{itemize}
\item $G_{\infty} \otimes_{\bC[\theta^{-1}]}\bC[\theta, \theta^{-1}] \simeq G $
\item $G_0 = \left( G_0 \cap G_{\infty} \right) \oplus \theta G_0$
\item $\tau\nabla_{\frac{\partial}{\partial \tau}}(G_{\infty}) \subset G_{\infty}$
\end{itemize}
(recall that $\tau=\theta^{-1}$).

\medskip
 
Even more concretely, and this is the way we will do it in Section \ref{Sec.: Frobenius for Quadric}, it can be done as follows. We need to find a $\mathbb{C}[\theta]$-basis of $G_0$ such that the connection matrix in this basis takes the form
\begin{align}\label{Eq.: Connection matrix Birkhoff I}
\left(\frac{A_0}{\theta}+A_{\infty}\right)\frac{d\theta}{\theta},
\end{align}
where $A_0$ and $A_{\infty}$ are constant matrices. To get the desired extension one just needs to consider these basis elements inside $G$ and define $G_{\infty}$ as a $\mathbb{C}[\theta^{-1}]$-submodule generated by them.

\subsubsection{Pairing.} 
\label{SubSubSec.: Pairing}

If $h$ is a cohomologically tame function, then there exists a non-degenerate bilinear pairing (cf. \cite{DoSa1}) 
\begin{align}\label{Eq.: Pairing II}
S_W \colon  G \otimes j^*G\to \bC[\theta, \theta^{-1}],
\end{align}
where $j\colon W \to W$ is given by $\theta \mapsto -\theta$.\footnote{The morphism $j$ extends uniquely to a morphism $\mathbb{P}^1_{\theta} \to \mathbb{P}^1_{\theta}$. Abusing notation we will denote this morphism again by~$j$, and we will use the same notation for its restrictions if it does not lead to confusion.} 

It satisfies
\begin{align}
\label{Eq.: Derivative Of Pairing}
&\frac{d}{d\theta}S_W(g_1, g_2)=S_W(\partial_{\theta}g_1,g_2)+S_W(g_1,\partial_{\theta} g_2),
\end{align}
i.e. it is a horizontal section of the sheaf $\mathcal{H}om_{\mathcal{O}_W}(\mathcal{F}_W \otimes j^*\mathcal{F}_W,\mathcal{O}_W) $ equipped with its natural connection, and

\begin{align}
\label{Eq.: Swap Arguments in the Pairing}
&S_W(g_1, g_2)=(-1)^n \overline{S_W(g_2, g_1)},
\end{align}
where we used the notation $\overline{P(\theta, \theta^{-1})}:=P(-\theta, -\theta^{-1})$ for a Laurent polynomial $P(\theta, \theta^{-1})$.

\medskip

Moreover, \eqref{Eq.: Pairing II} has the property 
\begin{align}\label{Eq.: teta^n}
S_W(G_0,j^*G_0) \subset \theta^n\bC[\theta] \subset \bC[\theta, \theta^{-1}],
\end{align}
and therefore we get a natural extension
\begin{align}\label{Eq.: Pairing III}
&S_{U_0} \colon  G_0 \otimes j^*G_0 \to \bC[\theta].
\end{align}

On $G_0$ we can write 
$$S_{U_0}=\sum_{i\geq n} S_i \theta^i,$$
where $S_i\colon G_0 \otimes j^*G_0 \to \bC \, \theta^i$ are \textit{higher residue pairings} of K. Saito; $S_n$ is the Grothendieck residue pairing. For a modern overview of K. Saito's works on this subject we refer to \cite{SaTa}.

\smallskip

Let $(\mathcal{F}, \nabla)$ be an extension as in Section \ref{SubSubSec.: Extension to P1}. We would like to extend (\ref{Eq.: Pairing III}) to a pairing
\begin{align}\notag
S \colon \mathcal{F} \otimes j^*\mathcal{F} \to \mathcal{O}_{\mathbb{P}^1_{\theta}}.
\end{align}

There exists $d\in \mathbb{Z}$ such that $S_W(G_{\infty},j^*G_{\infty}) \subset \tau^{-d} \mathcal{O}_{U_{\infty}}$ and therefore (\ref{Eq.: Pairing III}) always extends to 
\begin{align}\label{Eq.: Pairing IV}
S \colon  \mathcal{F} \otimes j^*\mathcal{F} \to \mathcal{O}_{\mathbb{P}^1_{\theta}}(- n \cdot \{0\} + d \cdot \{\infty \}).
\end{align}
Here $\mathcal{O}_{\mathbb{P}^1_{\theta}}(-n \cdot \{0\} + d \cdot \{\infty \})$ is an invertible subsheaf of $K_{\mathbb{P}^1_{\theta}}$ which consists of rational functions of ${\mathbb{P}^1_{\theta}}$, generated by $\theta^{n}$ and $\tau^{-d}$. It is isomorphic to $\mathcal{O}_{\mathbb{P}^1_{\theta}} $ if and only if $d=n$. The choice of $d$ in (\ref{Eq.: Pairing IV}) is not unique but there exists the minimal possible $d$. 

By (\ref{Eq.: teta^n}) we know that $d\geq n$ and therefore (\ref{Eq.: Pairing IV}) produces a pairing with values in $\mathcal{O}_{\mathbb{P}^1_{\theta}}$ iff $d=n$. The latter condition is equivalent to the existence of a global basis $e_1, \dots, e_{\mu}$ of $\mathcal{F}$ such that $S_{U_0}({e_i}_{|U_0}, {e_j}_{|U_0})\in \theta^n \mathbb{C}$.

\subsubsection{$V$-filtration.} 
\label{SubSubSec.: V-filtration}

Let $X$ be a smooth algebraic variety, $Y$ its closed smooth subvariety of codimension one, and $I$ the ideal sheaf of $Y$ in $X$. First define an increasing filtration $V_{\bullet} \mathcal{O}_X$ by putting $V_{i} \mathcal{O}_X=\mathcal{O}_X$ if $i\geq 0$ and $V_{i} \mathcal{O}_X=I^{-i}$ if $i < 0$. Now let $V_{\bullet} \mathcal{D}_X$ be an increasing filtration defined as
\begin{align*}
V_i\mathcal{D}_X=\{P\in \mathcal{D}_X \, | \, P(V_{m} \mathcal{O}_X)\subset V_{m+i} \mathcal{O}_X, \quad \forall \,\, m\in \mathbb{Z}\}.
\end{align*}
One can locally describe it more explicitly as follows (cf. \cite{PeSt}). Let $(y_1,\dots,$ $ y_n, x)$ be a local coordinate system on $X$ such that in this neighbourhood  $Y$ is given by the equation $x=0$. Then $V_0\mathcal{D}_X$ is a subsheaf of rings of $\mathcal{D}_X$ locally generated by $\mathcal{O}_X$, vector fields $\frac{\partial}{\partial y_1}, \dots , \frac{\partial}{\partial y_n}$ and $x\frac{\partial}{\partial x}$. If we denote $\partial_{x}=\frac{\partial}{\partial x}$, then $V_i\mathcal{D}_X$ is a $V_0\mathcal{D}_X$-module generated by $x^i\partial_{x}^j$ with $i-j\geq -k$.

\smallskip

Let $\mathcal{M}$ be a (left) $\mathcal{D}_X$-module and $V_{\bullet}\mathcal{M}$ a 
discrete exhaustive increasing filtration indexed by $\mathbb{Q}$. It is called \textit{$V$-filtration} iff

\smallskip

1.  it is compatible with the filtration $V_{\bullet}\mathcal{D}_X$, i.e. $\left(V_i\mathcal{D}_X\right)\left(V_{\alpha}\mathcal{M}\right) \subset V_{\alpha+i}\mathcal{M}$ for all $\alpha$ and~$i$; furthermore, the inclusion $I \left( V_{\alpha}\mathcal{M} \right) \subset V_{\alpha-1}\mathcal{M}$ should be an equality for $\alpha < 0$.

\smallskip

2. the action of $x \partial_{x}+\alpha$ on $\text{Gr}_{\alpha}^V\mathcal{M}$ is nilpotent.

\medskip

\noindent If such a filtration exists, then it is unique (cf. \cite{Bu}).

\medskip

The Gauss-Manin system $G$ considered as a $\mathbb{C}[\tau]\langle \partial_{\tau}\rangle$-module\footnote{This just means that we consider not $G$ itself but its push-forward as a $D$-module with respect to the open inclusion $U_0\cap U_{\infty} \to U_{\infty}$.} always has a $V$-filtration along $\{\tau=0\}$, and pairing (\ref{Eq.: Pairing II}) satisfies
\begin{align}\label{Eq.: Pairing Property wrt V-filtration}
S_W(V_0G, \overline{V_{<1}G}) \subset \mathbb{C}[\tau].
\end{align}
For more details we refer to \cite{DoSa1}.

\subsubsection{Tameness.} 
\label{SubSubSec.: Tameness}

Let $X$ be a smooth algebraic variety and $h \colon X \to \mathbb{A}^1$ a morphism. By a partial compactification we mean a commutative diagram
\begin{align}\notag
\xymatrix{
X \ar[r]^j \ar[d]^h &  \overline{X} \ar[dl]^{\overline{h}}\\
\mathbb{A}^1
}
\end{align}
where $\overline{X}$ is an algebraic variety(not necessarily smooth), $j$ is an open embedding, and $\overline{h}$ is proper.

\smallskip

The morphism $h$ is called \textit{cohomologically tame} iff there exists a partial compactification such that the support of $\Phi_{\overline{h}-a}(Rj_*\mathbb{C}_X)$ is finite and contained in $X_a$, for all $a \in \mathbb{A}^1$. We refer to \cite{Sa2} for more details.

\subsection{Initial conditions.}
\label{SubSec.: Initial conditions}

Let $(M,\, \circ,\, e,\, g,\, E)$ be a Frobenius manifold with an Euler field. In this setting one defines two endomorphisms of $\cT_M$ as
\begin{align}\label{Eq.: Operator mathcalUandV}
&\mathcal{U}(X)=E \circ X, \hspace{30pt} \mathcal{V}(X)=\nabla_X(E)-\frac{D}{2}X,
\end{align}
where $D$ is defined by  $Lie_E(g)=D g$ (see \cite[II.1]{Ma}).

If $p\in M$ is a semi-simple point of $M$, i.e. the algebra $(T_pM, \, \circ_p)$ is semi-simple (isomorphic to $\mathbb{C}^n$), then in a neighborhood of this point the tuple $(M,\, \circ,\, e,\, g,\, E)$ is uniquely determined by the data
\begin{align}\label{Eq.: Initial Conditions General Form}
(T, \, \cU, \, \cV, \, g, \, e),
\end{align}
where $T=T_pM$, $\cU$ and $\cV$ are endomorphisms of $T$ induced by \eqref{Eq.: Operator mathcalUandV}, $g$ is a non-degenerate symmetric bilinear pairing on $T$  induced by the metric, and $e$ is an element in $T$ induced by the identity vector field. This follows from \cite[Main Th., p.188]{Du} or \cite[Th. VII.4.2]{Sa1}.

\section{Construction of Landau-Ginzburg potentials}
\label{Sec.: LG Potentials}

We start by summarizing some facts about quantum cohomology of a smooth three-dimensional quadric~$Q_3$ (see~\cite{BaMa} for details). Then we explain how to obtain the standard LG potential for $Q_3$ from its quantum cohomology. As we already mentioned, this LG potential does not have enough critical points to be an honest LG model in the sense of Section \ref{SubSec.: LG models}, and we present its \textit{adhoc} partial compactification.

\subsection{Quantum cohomology of $Q_{3}$.} 

Let $V=Q_{3}$ be a smooth Fano hypersurface in $\mathbb{P}^{4}$ which is given by a non-degenerate homogeneous polynomial of degree 2. The singular cohomology groups $H^i(V, \mathbb{Z})$ are free of rank one in each even degree and vanish in odd degrees. Consider a graded basis $\Delta_0, \Delta_1, \Delta_2,\Delta_3$ of  $H^*(V,\mathbb{Z})$, such that $\Delta_0$ is the identity, $\Delta_1$ is the hyperplane class, $\Delta_1 \cup \Delta_{2}=\Delta_{3}$, where $\Delta_{3}$ is Poincar\'e dual to the class of a point. 

The table of quantum multiplication by $\Delta_1$ in the small quantum cohomology is
\begin{align}\label{Eq.: Multiplication Table Odd Quadric}
&\Delta_1\Delta_0=\Delta_1\\ \notag
&\Delta_1\Delta_{1}=2\Delta_{2}\\ \notag
&\Delta_1\Delta_{2}=\Delta_{3}+q\Delta_0\\ \notag
&\Delta_1\Delta_{3}=q\Delta_1.
\end{align}
Hence, the spectral cover consists of $4$ reduced points
\begin{align*}
&P_0=(1,0, 0, -q)\\ 
&P_i=(1, \xi_i,\frac{\xi_i^{2}}{2} , q ),
\end{align*}
where $\xi_i$ are roots of $\xi^{3}=4q$, and $1 \leq i \leq 3$. The point $P_0$ does not lie on the torus $\textbf{T}$ (cf. Section \ref{SubSubSec.: Spectral Cover for QH}).

\subsubsection{Initial conditions.}
\label{SubSubSec.: InitCondQuad}

We can express the anti-canonical class as
$$
-K_V=3\D_1.
$$
In the basis of the $\D_i$'s the initial conditions take the form
\begin{align*}
\cU=\left(
  \begin{array}{rrrr}
    0 & 0 & 3q & 0 \\
    3 & 0 & 0  & 3q \\
    0 & 6 & 0  & 0 \\
    0 & 0 & 3 & 0 \\
  \end{array}
\right)\\
\intertext{and}
\cV=\left(
  \begin{array}{rrrr}
    1 & 0 & 0 & 0 \\
    0 & 0 & 0  & 0 \\
    0 & 0 & -1  & 0 \\
    0 & 0 & 0 & -2 \\
  \end{array}
\right)+\frac{1}{2}.
\end{align*}

\subsection{Standard LG potential.} 

Restricting to the torus $\textbf{T}$ we can rewrite system \eqref{Eq.: Multiplication Table Odd Quadric} as
\begin{align*}
&\Delta_1=\frac{2\Delta_2}{\Delta_1}\\ 
&\Delta_1=\frac{\Delta_{3}+q}{\Delta_{2}}\\ 
&\Delta_1=\frac{q\Delta_1}{\Delta_{3}}.
\end{align*}
The above system can be rewritten as 
\begin{align}\label{Eq.: Multiplication table III}
&\Delta_1=\frac{2\Delta_2}{\Delta_1}\\ \notag
&\frac{\Delta_{2}}{\Delta_{1}}=\frac{(\Delta_{3}+q)^2}{2\Delta_{2}\Delta_{3}} \\ \notag
&\frac{\Delta_{3}^2}{2\Delta_{2}\Delta_{3}}=\frac{q^2}{2\Delta_{2}\Delta_{3}}.
\end{align}
It is easy to see that if we define
\begin{align}\label{Eq.: Non-compactified potential I}
f=\Delta_1+\frac{2\Delta_2}{\Delta_1} + \frac{(\Delta_{3} + q)^2}{2\Delta_{2}\Delta_{3}},
\end{align}
then the system $\D_i\frac{\partial f}{\partial \D_i}=0$ coincides with \eqref{Eq.: Multiplication table III}. In this sense $f$ "integrates" the multiplication table.

\medskip

We claim that \eqref{Eq.: Non-compactified potential I} is the standard LG potential for $Q_3$ proposed in \cite{EgHoXi}. Indeed, consider another coordinate system on the torus $\textbf{T}=\{\D_1\D_2 \D_{3} \neq 0\}$ given by 
$$
Y_1=\Delta_1, \quad Y_2=\frac{2\Delta_2}{\Delta_1}, \quad Y_{3}=\Delta_{3}.
$$
Rewriting \eqref{Eq.: Non-compactified potential I} in terms of these coordinates we get
\begin{align}\label{Eq.: Non-compactified potential II}
f=Y_1+Y_{2}+\frac{(Y_{3}+q)^2}{Y_1Y_2Y_3}.
\end{align}
which is exactly the LG potential proposed in loc.cit..

\subsection{Compactification.} 

By construction \eqref{Eq.: Non-compactified potential I} has 3 critical points but a Landau-Ginzburg potential for $Q_3$ in the sense of Section \ref{SubSec.: LG models} must necessarily have 4 critical points. Below we will give an \textit{adhoc} partial compactification of $f$ to a new LG potential $\widetilde{f}$ which has the correct number of critical points. In Section \ref{Sec.: Frobenius for Quadric} we will study the Gauss-Manin system of $\widetilde{f}$ and show that it deserves the name LG potential.

\smallskip

Consider the affine space $\mathbb{A}^3=\Spec\bC[Y_1, Y_2, Y_3]$. Expression \eqref{Eq.: Non-compactified potential II} gives a regular function on the torus $\{Y_1Y_2Y_3 \neq 0\} \subset \bA^3$. Functions $x, y, z$ given by
\begin{align} \label{Eq.: Coordinate Change}
&x=\frac{Y_3+q}{qY_1}\\  \notag
&y=Y_1\\ \notag
&z=\frac{Y_2}{Y_1}-1
\end{align}
define another coordinate system on this torus. Rewriting $f$ in terms of these coordinates we get
\begin{align}\label{Eq.: Compactified LG model without parameters}
&\widetilde{f}=y(2+z)+\frac{qx^2}{(xy-1)(1+z)}.
\end{align}
One can interpret \eqref{Eq.: Compactified LG model without parameters} as a regular function on an open subvariety of $\mathbb{A}^{3}=\Spec\bC[x, y,z]$ defined by $\{(xy-1)(1+z)\neq 0\}$. The torus $\{Y_1Y_2Y_{3} \neq 0\}$ is embedded into this space by formulas \eqref{Eq.: Coordinate Change}. 

It is easy to check that the critical locus of $\widetilde{f}$ consists of $4$ points
$$
P_0=(0, 0, -2) \quad \text{and} \quad P_i=(\frac{2}{\xi_i},\xi_i,0),
$$
where $\xi_i^3=4q$.

\subsection{General case.}

The above considerations work for an arbitrary smooth odd-dimensional quadric $Q_{2n+1}$ of dimension $2n+1$. Let $\Delta_0, \dots, \Delta_{2n+1}$ be a graded basis of  $H^*(Q_{2n+1},\mathbb{Z})$, such that $\Delta_0$ is the identity, $\Delta_1$ is the hyperplane class, $\Delta_i=\Delta_1^{\cup i}$ for $i\leq n$, and $\Delta_i \cup \Delta_{2n+1-i}=\Delta_{2n+1}$, where $\Delta_{2n+1}$ is Poincar\'e dual to the class of a point. 

One can write down the quantum multiplication by $\D_1$. Then one can show that 
\begin{align*}
f=\Delta_1+\frac{\Delta_2}{\Delta_1}+\dots +\frac{\Delta_n}{\Delta_{n-1}}+\frac{2\Delta_{n+1}}{\Delta_n}+\frac{\Delta_{n+2}}{\Delta_{n+1}} +\dots + \frac{\Delta_{2n}}{\Delta_{2n-1}} + \frac{(\Delta_{2n+1} + q)^2}{2\Delta_{2n}\Delta_{2n+1}}
\end{align*}
is the analogue of \eqref{Eq.: Non-compactified potential I}, and
\begin{align*}
&\widetilde{f}=\sum_{i=1}^ny_i(2+z_i)+\frac{qx^2}{(xy_1\dots y_n-1)(1+z_1)\dots (1+z_n)}.
\end{align*}
is the analogue of \eqref{Eq.: Compactified LG model without parameters}.

\section{Gauss-Manin system of $\widetilde{f}$}
\label{Sec.: Frobenius for Quadric}

\subsection{Notation.}

For convenience we repeat the partial compactification in a somewhat backwards order. Namely, consider $\bA^3=\Spec \bC[x,y,z]$ and let $\widetilde{U}$ be the open subvariety defined by $\{(xy-1)(1+z)\neq 0\}$. On $\widetilde{U}$ we have the regular function 
\begin{align}\label{Eq.: F for Q_3}
&\widetilde{f}=y(2+z)+\frac{qx^2}{(xy-1)(1+z)},
\end{align}
which is our partially compactified potential (\ref{Eq.: Compactified LG model without parameters}).

Consider functions $\D_1, \D_2, \D_3$ on $\bA^3$ given by
\begin{align}\label{Eq.: Formulas for Delta_i}
\D_1=y, \quad \D_2=\frac{y^2}{2}(z+1), \quad \D_3=qxy-q,
\end{align}
which form a coordinate system on the subset $\{y \neq 0\} \subset \bA^3$. The inverse coordinate change is given by 
\begin{align}\label{Eq.: xyz via deltas}
x=\frac{\D_3+q}{q\D_1}, \quad y=\D_1, \quad z=\frac{2\D_2}{\D_1^2}-1.
\end{align}

Let $U$ be the intersection $\widetilde{U} \cap \{y \neq 0\}$. On $U$ function (\ref{Eq.: F for Q_3}) can be rewritten in terms of $\D_1, \D_2, \D_3$ as 
\begin{align}\label{Eq.: f for Q_3}
f:=\widetilde{f}_{|U}=\D_1+\frac{2\D_2}{\D_1}+\frac{(\D_3+q)^2}{2\D_2\D_3}.
\end{align}
Formulas \eqref{Eq.: Formulas for Delta_i} give an isomorphism of $U$ with the algebraic torus $\textbf{G}_m^3=\Spec \bC[t_1^{\pm 1}, t_2^{\pm 1}, t_3^{\pm 1}]$, such that $t_i$'s correspond to $\D_i$'s. Formula \eqref{Eq.: f for Q_3} gives the LG potential before the compactification as in \eqref{Eq.: Non-compactified potential I}.

\subsubsection{Lemma.}
\label{SubSubSec.: Lemma Tameness}

{\it Function (\ref{Eq.: F for Q_3}) is cohomologically tame.\footnote{It is also true that $f$ has isolated singularities at infinity in the sense of \cite{Do2}.}
}

\smallskip

\proof See Appendix A. $\blacksquare$

\subsubsection{Lemma.} 

{\it The Gauss-Manin system  $G^{\widetilde{f}}$ has the following properties:

\smallskip

(i) $G^{\widetilde{f}}$ is a free $\mathbb{C}[\theta,\theta^{-1}]$-module of rank 4;

\smallskip

(ii) $G_0^{\widetilde{f}}$ is a free $\mathbb{C}[\theta]$-module of rank 4.
}

\smallskip

\proof For both properties it is essential that $\widetilde{f}$ is cohomologically tame.

\smallskip

(i) For a function with isolated critical points the module $G$ is always free of finite rank. If, moreover, the function is cohomologically tame, then the rank is equal to the Milnor number (\cite{Do2}, Th. 5.2.3). In our case it is 4.

\smallskip

(ii) For a function with (cohomologically) isolated critical points at infinity Corollary 5.2.6 of \cite{Do2} states, that $G_0$ is free and of finite type iff the function is cohomologically tame. 

Applying this corollary to $\widetilde{f}$ we get that $G_0^{\widetilde{f}}$ is a free $\mathbb{C}[\theta]$-module of finite rank. Hence, its rank equals to the dimension of the fiber at zero. Using Proposition 5.1.1 of \cite{Do2} we see that the rank is equal to the Milnor number.~$\blacksquare$

\subsubsection{Lemma.} 
\label{SubSubSec.: GM Iso}

{\it The natural morphism of $\mathcal{D}_W$-modules $G^{\widetilde{f}} \to G^f$ given by the restriction of differential forms from $\widetilde{U}$ to $U$ is an isomorphism.\footnote{Recall from Section \ref{SubSec.: GM system} that $W=\mathbb{A}^1_{\theta}-\{0\}$}
}

\smallskip

\proof Restriction of differential forms from $\widetilde{U}$ to $U$ defines the morphism 
$$
\Omega^i(\widetilde{U}) \to \Omega^i(U),
$$
which is injective but not surjective; it is the localization morphism given by inverting $\D_1$. One can check directly that the induced morphism $G^{\widetilde{f}} \to G^f$ on the Gauss-Manin systems is also injective.

\smallskip

By Theorem 5.2.3 of \cite{Do2} the rank of $G^f$ is 4 (we use here that $f$ has one isolated singularity at infinity).

\smallskip

Consider the short exact sequence of $\mathcal{O}_W$-coherent $\mathcal{D}_W$-modules 
\begin{align*}
0 \to G^{\widetilde{f}} \to G^f \to G^f/G^{\widetilde{f}} \to 0.
\end{align*}
Since $\text{rk }G^{\widetilde{f}}=\text{rk }G^f$ the quotient is an $\mathcal{O}_W$-module of rank zero. Therefore, by the standard fact that for a smooth algebraic variety $X$ any $\mathcal{O}_X$-coherent $\mathcal{D}_X$-module is a locally free $\mathcal{O}_X$-module (see \cite{Be}, Lect. 2, 1.a), we get that $G^f/G^{\widetilde{f}}$ is locally free of rank zero and hence vanishes. $\blacksquare$

\subsection{Birkhoff problem.} 
\label{SubSec.: Our Birkhoff problem }

Consider the following 3-form on $\widetilde{U}$
\begin{align}\notag
&\omega_0=\frac{dx\wedge dy \wedge dz}{(xy-1)(z+1)},
\end{align}
and let $\omega_i=\D_i \omega_0$. Note also that
\begin{align}\notag
&{\omega_0}_{|U}=\frac{d \D_1}{\D_1}\wedge\frac{d \D_2}{\D_2}\wedge\frac{d \D_3}{\D_3}.
\end{align}
If $\omega$ is a 3-form, then let $[\omega]$ denote its class in $G_0$. In the above formulas by $\D_i$ we mean ${\D_i}_{|\widetilde{U}}$ and ${\D_i}_{|U}$ respectively. We will continue to use this notation, if it does not lead to confusion.

\subsubsection{Lemma.}
\label{SubSubSec.: Identities in G} 

{\it In $G^f$ we have the following identities}
\begin{align*}
&[\D_if'_{\D_i}\omega_0]=0\\
&[\D_i\D_jf'_{\D_i}\omega_0]=0\\
&[\D_i^2f'_{\D_i}\omega_0]=\theta[\omega_i].
\end{align*}

\proof Let us only prove the third identity for $i=2$. The other cases are analogous. 

We have the following equality of differential forms
\begin{align*}
\D_2^2f'_{\D_2}\omega_0=-df\wedge \left(\D_2\frac{d \D_1}{\D_1}\wedge\frac{d \D_3}{\D_3}\right),
\end{align*}
hence $[\D_2^2f'_{\D_2}\omega_0]=\theta[\omega_2]$ in $G^f$. $\blacksquare$

\subsubsection{Lemma.} 
\label{SubSubSec.: Lemma Linear Independence}

{\it Elements $[\omega_0], \dots, [\omega_3]$ are $\mathbb{C}[\theta]$-linearly independent in $G_0^{\widetilde{f}}$.}

\smallskip

\proof The vector space $G_0^{\widetilde{f}}/\theta G_0^{\widetilde{f}}$ can be identified with the Milnor ring by mapping $1$ to the class of $[\omega_0]$. Under this isomorphism the class of $\D_i$ goes to the class of $[\omega_i]$. Since $1, \D_1, \D_2, \D_3$ form a basis in the Milnor ring, classes of $[\omega_0], \dots, [\omega_3]$ form a basis in $G_0^{\widetilde{f}}/\theta G_0^{\widetilde{f}}$. This implies the statement. $\blacksquare$

\subsubsection{Lemma.} 
\label{SubSubSec.: Lemma A_0, A_inf}

{\it
(i) Elements $[\omega_0], \dots, [\omega_3]$ freely generate in $G^{\widetilde{f}}$ an $\mathcal{O}_W$-submodule $H^{\widetilde{f}}$ of rank 4;

\smallskip

(ii) The following identities hold
\begin{align*}
&\theta^2 \partial_{\theta}[\omega_0]=3[\omega_1]\\
&\theta^2 \partial_{\theta}[\omega_1]=6[\omega_2]+\theta[\omega_1]\\
&\theta^2 \partial_{\theta}[\omega_2]=3[\omega_3]+3q[\omega_0]+2\theta[\omega_2]\\
&\theta^2 \partial_{\theta}[\omega_3]=3q[\omega_1]+3\theta[\omega_3],
\end{align*}
and therefore $H^{\widetilde{f}}$ is a $\mathcal{D}_W$-submodule; 

\smallskip

(iii) $G^{\widetilde{f}}=H^{\widetilde{f}}$.

\smallskip

(iv) The connection matrix in the basis $[\omega_0], \dots, [\omega_3]$ takes the form
\begin{align}\notag
\left(\frac{A_0}{\theta}+A_{\infty}\right)\frac{d\theta}{\theta},
\end{align}
where
\begin{align}\notag
A_0=
\left(
  \begin{array}{rrrr}
    0 & 0 & 3q & 0 \\
    3 & 0 & 0  & 3q \\
    0 & 6 & 0  & 0 \\
    0 & 0 & 3 & 0 \\
  \end{array}
\right)
\intertext{and}\notag
A_{\infty}=
\left(
  \begin{array}{rrrr}
    0 & 0 & 0 & 0 \\
    0 & 1 & 0  & 0 \\
    0 & 0 & 2  & 0 \\
    0 & 0 & 0 & 3 \\
  \end{array}
\right)
\end{align}
}

\medskip

\proof (i) By Lemma \ref{SubSubSec.: Lemma Linear Independence} $[\omega_0], \dots, [\omega_3]$ are linearly independent in 
$G_0^{\widetilde{f}}$, and hence also in $G^{\widetilde{f}}$ and $G^f$. Therefore, they generate a submodule of rank 4 in $G^{\widetilde{f}}$ (and in $G^f$).

\smallskip

(ii) Because of the natural isomorphism $G^{\widetilde{f}} \to G^f$ we can check these identities in $G^f$.

First, note that the following identities hold in the ring of functions on $U$
\begin{align*}
&f=3\D_1-2\D_1f'_{\D_1}-\D_2f'_{\D_2}\\
&\D_1\D_1=2\D_2+\D_1^2f'_{\D_1}\\ 
&\D_1\D_2=(\D_3+q)+\D_2(\D_1f'_{\D_1}+\D_2f'_{\D_2}-\D_3f'_{\D_3})\\
&\D_1\D_3=q\D_1+\D_3(\D_1f'_{\D_1}+\D_2f'_{\D_2}+\D_3f'_{\D_3})-q(\D_1f'_{\D_1}+\D_2f'_{\D_2}-\D_3f'_{\D_3}).
\end{align*}
These identities can be checked by direct computations.

Using the first identity we get
\begin{align*}
&\theta^2 \partial_{\theta}[\omega_0]=[f\omega_0]=[(3\D_1-2\D_1f'_{\D_1}-\D_2f'_{\D_2})\omega_0]=\\
&3[\D_1\omega_0]-2[\D_1f'_{\D_1}\omega_0]-[\D_2f'_{\D_2}\omega_0].
\end{align*}
Applying Lemma \ref{SubSubSec.: Identities in G} we get
\begin{align*}
&\theta^2 \partial_{\theta}[\omega_0]=3[\D_1\omega_0]=3[\omega_1].
\end{align*}

Using the first two identities and Lemma \ref{SubSubSec.: Identities in G} we get 
\begin{align*}
&\theta^2 \partial_{\theta}[\omega_1]=[f\D_1\omega_0]=[6\D_2\omega_0+\D_1^2f'_{\D_1}\omega_0-\D_1\D_2f'_{\D_2}\omega_0]=6[\omega_2]+\theta[\omega_1].
\end{align*}

The remaining two formulas are obtained analogously.

\smallskip

(iii) Since $H^{\widetilde{f}}$ and $G^{\widetilde{f}}$ are $\mathcal{O}_W$-coherent $\mathcal{D}_W$-modules of the same rank, they coincide (as in the proof of Lemma \ref{SubSubSec.: GM Iso}). 

\smallskip

(iv) It follows from (ii). $\blacksquare$

\subsubsection{Lemma.} 
\label{SubSubSec.: Lemma about Basis}

{\it The classes $[\omega_0],\dots, [\omega_3]$ form a $\mathbb{C}[\theta]$-basis in $G_0^{\widetilde{f}}$.}

\medskip

\proof Let $H_0^{\widetilde{f}}$ be the $\mathcal{O}_{\mathbb{A}^1_{\theta}}$-submodule of $ G_0^{\widetilde{f}}$ generated by $[\omega_0], \dots, [\omega_3]$. We have the short exact sequence of $\mathcal{O}_{\mathbb{A}^1_{\theta}}$-modules
\begin{align}\label{Eq.: ShortExactSequence 0}
0 \to H_0^{\widetilde{f}} \to G_0^{\widetilde{f}} \to Q_0^{\widetilde{f}} \to 0,
\end{align}
and we need to show that $Q_0^{\widetilde{f}}=0$.

Since ${Q_0^{\widetilde{f}}}_{|\mathbb{A}^1_{\theta}-\{0\}}=0$ by Lemma \ref{SubSubSec.: Lemma A_0, A_inf}, and $Q_0^{\widetilde{f}}$ is finitely generated, it is enough to prove that the fiber at zero vanishes, i.e. $Q_0^{\widetilde{f}}\otimes_{\mathbb{C}[\theta]} \mathbb{C}[\theta]/(\theta)=0$. 

\smallskip

Tensoring \eqref{Eq.: ShortExactSequence 0} with $\mathbb{C}[\theta]/(\theta)$ we get a short exact sequence (tensor product is right exact)
\begin{align*}
H_0^{\widetilde{f}}\otimes_{\mathbb{C}[\theta]} \mathbb{C}[\theta]/(\theta) \to G_0^{\widetilde{f}} \otimes_{\mathbb{C}[\theta]} \mathbb{C}[\theta]/(\theta)\to Q_0^{\widetilde{f}}\otimes_{\mathbb{C}[\theta]} \mathbb{C}[\theta]/(\theta) \to 0,  
\end{align*}
which can be rewritten as 
\begin{align*}
H_0^{\widetilde{f}}\otimes_{\mathbb{C}[\theta]} \mathbb{C}[\theta]/(\theta) \to \Omega^n(\widetilde{U})/d\widetilde{f}\wedge \Omega^{n-1}(\widetilde{U}) \to Q_0^{\widetilde{f}}\otimes_{\mathbb{C}[\theta]} \mathbb{C}[\theta]/(\theta) \to 0.
\end{align*}
Since classes of $[\omega_0], \dots, [\omega_3]$ generate $\Omega^n(\widetilde{U})/d\widetilde{f}\wedge \Omega^{n-1}(\widetilde{U})$, the first map is surjective. Therefore,  $Q_0^{\widetilde{f}}\otimes_{\mathbb{C}[\theta]} \mathbb{C}[\theta]/(\theta)=0$, and, finally,  $Q_0^{\widetilde{f}}=0$. $\blacksquare$

\subsection{Pairing.} 
\label{SubSec.: Our pairing}

In this section we study pairing \eqref{Eq.: Pairing II} in our setup. Since it will make no difference here, we are dropping the subscripts in the notation of the Gauss-Manin systems, and just write $G$.

\subsubsection{Lemma.} 

{\it The $V$-filtration on $G$ along $\{\tau=0\}$ is given by
\begin{align*}
&V_0G=\bigoplus_{i=0}^{3}\mathbb{C}[\tau]e_i\\ 
&V_pG=\tau^{-p}V_0G,
\end{align*}
where $e_i=\tau^i[\omega_i]$.
}

\medskip

\proof 
This lemma, as well as most of the results of this article, first appeared in \cite{Sm}. There the proof of this lemma has a mistake which is corrected here.

\smallskip

It is enough to show that this filtration satisfies the conditions of Section \ref{SubSubSec.: V-filtration}.

\smallskip

\textit{1. Compatibility of filtrations.}

\textit{1a.} It is clear that $\tau (V_pG) \subset V_{p-1}G$, and using that $\partial_{\tau}=-\theta^2 \partial_{\theta}$ and applying Lemma \ref{SubSubSec.: Lemma A_0, A_inf} it is not difficult to see that $\partial_{\tau} (V_pG) \subset V_{p+1}G$.

\smallskip

These two facts imply that $(V_m\mathcal{D}_{\mathbb{A}_{\tau}})(V_pG) \subset V_{p+m}G$.

\smallskip

\textit{1b.} It is clear that the condition $\tau \, V_pG = V_{p-1}G$ for $p<0$ holds. 

\smallskip

\textit{2. Nilpotence.} Classes of $\tau^{-p}e_0,\dots,  \tau^{-p}e_3$ form a basis in $\text{Gr}_p^VG$. Using  Lemma \ref{SubSubSec.: Lemma A_0, A_inf} one can see that in this basis the operator induced by $\tau \partial_{\tau}+p$ on $\text{Gr}_p^VG$ is given by the matrix
\begin{align}\notag
\left(
  \begin{array}{rrrr}
    0  & 0  &  0 & 0 \\
    -3 & 0  &  0 & 0 \\
    0  & -6 &  0 & 0 \\
    0  & 0  & -3 & 0 \\
  \end{array}
\right)
\end{align}
It is clearly nilpotent. $\blacksquare$

\subsubsection{Lemma.} 
\label{SubSubSec.: PairingComputation}

{\it The pairing $S_W$ satisfies\footnote{The overline over the second argument of $S_W$ stresses that the element is considered as a section of $j^*\mathcal{F}_W$. Therefore, $\tau$ and $\nabla_{\frac{\partial}{\partial\tau}}$ act with the opposite sign.}
\begin{align}
S_W([\omega_k],\overline{[\omega_l]})=
\lbrace \begin{array}{lll}
S_W([\omega_0],\overline{[\omega_3]}) &  \text{if} \quad  k+l=3 \\
0                                   & \text{otherwise}  
\end{array}
\end{align}
and $S_W([\omega_0],\overline{[\omega_3]}) \in \tau^{-3} \mathbb{C}$.
}

\medskip

\proof To simplify the notation we will be writing $S$ instead of $S_W$. By \eqref{Eq.: teta^n} we know that 
$$
S([\omega_k],\overline{[\omega_l]})\in \tau^{-3}\mathbb{C}[\tau^{-1}].
$$
On the other hand, by \eqref{Eq.: Pairing Property wrt V-filtration} we get
$$
S(e_k,\overline{e_l})=\tau^k(-\tau)^l S([\omega_k],\overline{[\omega_l]})\in \mathbb{C}[\tau],
$$
therefore $S([\omega_k],\overline{[\omega_l]})\in \tau^{-(k+l)}\mathbb{C}[\tau]$. Hence 
\begin{align}\label{Eq.: Vanishing By V-filtration}
& S([\omega_k],\overline{[\omega_l]})=0 \hspace{37pt} \text{if} \quad k+l < 3, \\ \notag 
& S([\omega_k],\overline{[\omega_l]}) \in \tau^{-3} \mathbb{C} \hspace{20pt} \text{if} \quad k+l = 3. 
\end{align}
  
To show vanishing in the remaining 4 cases with $k+l > 3$ one just combines \eqref{Eq.: Derivative Of Pairing} and \eqref{Eq.: Vanishing By V-filtration}. Let us consider the case $k=1$, $l=3$ only. Applying \eqref{Eq.: Derivative Of Pairing} to $S([\omega_0], \overline{[\omega_3]})\in \tau^{-3}\mathbb{C}$  and using \eqref{Eq.: Vanishing By V-filtration} we get
\begin{align}\notag
&-3\tau^{-1}S([\omega_0], \overline{[\omega_3]})=\frac{d}{d\tau}S([\omega_0], \overline{[\omega_3]})=S(\partial_{\tau} [\omega_0], \overline{[\omega_3]})-S([\omega_0], \overline{\partial_{\tau}[\omega_3]})=\\ \notag
& \hspace{85pt} =-S(3[\omega_1], \overline{[\omega_3]})+S([\omega_0], \overline{3q[\omega_1]+3\theta[\omega_3]})=\\ \notag
& \hspace{85pt} =-3S([\omega_1], \overline{[\omega_3]})-3\tau^{-1}S([\omega_0], \overline{[\omega_3]}),
\end{align}
and therefore $S([\omega_1], \overline{[\omega_3]})=0$.

\smallskip

Similarly one can show that $S([\omega_0], \overline{[\omega_3]})=S([\omega_2], \overline{[\omega_1]})$. Moreover, by \eqref{Eq.: Swap Arguments in the Pairing} we have $S([\omega_0], \overline{[\omega_3]})=S([\omega_3], \overline{[\omega_0]})$ and $S([\omega_2], \overline{[\omega_1]})=S([\omega_1], \overline{[\omega_2]})$. $\blacksquare$

\subsection{Canonical solution to the Birkhoff problem.}

The problem of extending $(G_0, \nabla)$ to $\mathbb{P}^1_{\theta}$ described in Section \ref{SubSec.: GM system} has a canonical solution given by Hodge theory. Here we will show that our solution given by the basis $\omega_0, \dots, \omega_3$  is canonical in this sense. We more or less keep notation of \cite[Sec.~5]{DoSa2}. Details can be found in loc. cit. and references therein.

It is a general fact (see \cite{Sa2}) that the vector space 
\begin{align}
H= \bigoplus_{\alpha \in [0,1) } \text{Gr}_{\alpha}^V G
\end{align}
carries a mixed Hodge structure, i.e. $H$ has a rational structure $H_{\bQ}$, an increasing weight filtration $W_{\bullet}H_{\bQ}$, a decreasing Hodge filtration $F^{\bullet}H$, s.t. the Hodge filtration induces a pure Hodge structure of weight $m$ on $\text{Gr}_m^WH$ for all $m$.

For the function $\wt{f}$ we have
\begin{align*}
&H=\text{Gr}_0^VG = \oplus_{i=0}^{3} \bC \, e_i \\ 
&F^{p}H=\bigoplus_{i=0}^{3-p} \bC \, e_i,
\end{align*}
where abusing notation we write $e_i$'s meaning classes of $e_i$'s in $H$. The complexification of the weight filtration is
\begin{align}\label{Eq.: Weight Filtration}
0 \subset \bC e_3 \subset \bC e_3 \subset \bC e_2 \oplus \bC e_3 \subset \bC e_2 \oplus \bC e_3 \subset \bC e_1 \oplus \bC e_2 \oplus \bC e_3 \subset \bC e_1 \oplus \bC e_2 \oplus \bC e_3 \subset H,
\end{align}
where the first term on the left is $W_{-1}H$. The only non-trivial associated graded objects are $\text{Gr}_0^WH$, $\text{Gr}_2^WH$ and $\text{Gr}_4^WH$.

In the case of the function $\wt{f}$ Saito's canonical opposite filtration (it is a filtration on $H$) is defined as
\begin{align*}
H_{\text{Saito}}^{\bullet}=\sum_q \overline{F^qH} \cap W_{3+q-\bullet}H.
\end{align*} 

From \eqref{Eq.: Weight Filtration} and the fact that the weight filtration is stable under conjugation we have
\begin{align}\label{Eq.: Conjugation}
\ol{e}_i=\sum_{r=i}^3 a_{ri}e_r,
\end{align}
with $a_{ii} \neq 0$. Using \eqref{Eq.: Conjugation} one can show that 
\begin{align}\label{Eq.: Saito's filtration}
H_{\text{Saito}}^{p}=\bigoplus_{i=0}^{3-p} \bC e_{3-i}.
\end{align}

To any solution of the Birkhoff problem one can attach a filtration $H^{\bullet}$ on $H$ by the formula
\begin{align*}
&H^i:=G_{\infty}^i \cap V_0G/G_{\infty}^i \cap V_{-1}G,
\end{align*}
where $G^{k}_{\infty}=\tau^k G_{\infty}$. To prove that this solution of the Birkhoff problem is canonical one needs to show that the filtrations $H^{\bullet}$ and $H_{\text{Saito}}^{\bullet}$ coincide.

In our case it is easy to see that
\begin{align*}
&H^i=\bigoplus_{i=0}^{3-p} \bC e_{3-i},
\end{align*}
which coincides with \eqref{Eq.: Saito's filtration}. Hence, our solution to the Birkhoff problem is canonical.

\subsection{Frobenius manifold.}

Ideally one would like to show the existence of (or to exhibit) a deformation of $\wt{f}$ producing a Saito's framework isomorphic to $QH(Q_3)$. We have not been able to achieve this goal so far. 

There is a general construction of such Saito's frameworks due to A. Douai and C. Sabbah (see \cite{DoSa1}) but it requires some additional properties of $\wt{f}$ that we have not been able to check. Namely, one needs to show that $\wt{f}$ is $M$-tame. We refer to loc.cit. for details.

Assume that such Saito's framework $(M, \circ, e, g_{\omega}, E)$ exists. Then the initial conditions for $M$ at the origin are
$$
T=T_{x_0}M, \quad  \cU=A_0, \quad \cV=-A_{\infty}+\frac{3}{2}\text{Id}, \quad g_{\omega}, \quad e.
$$
Therefore, by Lemmas \ref{SubSubSec.: Lemma A_0, A_inf} and \ref{SubSubSec.: PairingComputation} we see that these initial conditions coincide with those of Section \ref{SubSubSec.: InitCondQuad}.

\appendix

\addtocontents{toc}{\protect\setcounter{tocdepth}{1}}

\section{Proof of Lemma \ref{SubSubSec.: Lemma Tameness}}

The proof of Lemma 4.1.1 given here has been kindly explained to us by Claude Sabbah and Andr\'as N\'emethi.

\subsection{Vanishing cycles.}

Here we recall some basic facts about functors of vanishing cycles. For a comlpex algebraic variety $X$ we denote by $D^b_c(X)$ the bounded derived category of $\mathbb{C}_{X^{an}}$-modules with constructible cohomology, where $X^{an}$ is the associated analytic space.

\subsubsection{Functor of vanishing cycles.} 

Let $X$ be a complex algebraic variety and $g \colon X \to \mathbb{A}^1_t$ a morphism. The functor of vanishing cycles to the fiber over $0$ of the morphism $g$ is denoted $\Phi_g$. If one considers the fiber over $a$, then one shifts $g$ and considers $\Phi_{g-a}$.

If we denote $X_0$ the fiber of $g$ over $0$, then the functor of vanishing cycles to this fiber is a triangulated functor
$$
\Phi_g \colon D^b_c(X) \to D^b_c(X_0),
$$
i.e. it maps distinguished triangles to distinguished triangles. See \cite{Di} for a precise definition.

\medskip

Below we collect some basic properties of these functors. They are standard and can be found in \cite{Di}.

\subsubsection{Proper morphism.} 

Let $\pi \colon X \to Y$ be a morphism of algebraic varieties and consider the following commutative diagram
\begin{align}\notag
\xymatrix{
X_0  \ar[r] \ar[d]^{\widehat{\pi}}    &   X  \ar[d]^{\pi}        \\
Y_0  \ar[r] \ar[d]          &   Y  \ar[d]^g          \\
\{0\}   \ar[r]                 &   \mathbb{A}^1
}
\end{align}
where $X_0$ and $Y_0$ are fibers over $0$. Naturally one can attach to it the diagram of derived categories and functors between them
\begin{align}\notag
\xymatrix{
D_c^b(X_0)   \ar[d]^{R\widehat{\pi}_*}    &  \ar[l]_{\Phi_{g \circ \pi}} D_c^b(X)  \ar[d]^{R\pi_*}        \\
D_c^b(Y_0)                              & \ar[l]_{\Phi_{g}}  D_c^b(Y)
}
\end{align}

\smallskip

{\bf  Fact:} if $\pi$ is \textit{proper}, then the above diagram is commutative (e.g. this is true if $\pi$ is a closed embedding).

\subsubsection{Restriction to an open subset.} 
\label{SubSubSec.: RestrictionToOpen}

Let $U \subset X$ be an open subset and let $j$ be the natural inclusion. We have a commutative diagram
\begin{align}\notag
\xymatrix{
U_0  \ar[r] \ar[d]^{\widehat{j}}    &   U  \ar[d]^{j}        \\
X_0  \ar[r] \ar[d]          &   X  \ar[d]^g          \\
\{0\}   \ar[r]                 &   \mathbb{A}^1
}
\end{align}
with $U_0$ and $X_0$ being fibers over $0$. Consider the associated diagram of derived categories and functors between them
\begin{align}\notag
\xymatrix{
D_c^b(U_0)       &  \ar[l]_{\Phi_{g \circ j}} D_c^b(U)        \\
D_c^b(X_0)  \ar[u]^{\widehat{j}^{-1}}                            & \ar[l]_{\Phi_{g}}  D_c^b(X) \ar[u]_{j^{-1}}
}
\end{align}

\smallskip

\textbf{Fact:} the above diagram is commutative.

\subsubsection{Duality.} 
\label{SubSubSec.: Duality For Appendix}

For any complex algebraic variety $Y$ there exists a functor $\textbf{D}_Y\colon D^b(Y) \to D^b(Y)$ defined by 
$$
\textbf{D}_Y(\mathcal{F}^{\bullet}):=R\mathcal{H}om(\mathcal{F}^{\bullet}, \omega_Y),
$$
where $\omega_Y$ is the dualizing complex. It has the property $\textbf{D}_Y\circ \textbf{D}_Y \simeq \text{Id}_{ D^b(Y)}$. Moreover, this functor restricts to $D^b_c(Y)$ and we will use the same notation for this restriction.  If we do not want specify the space, we will just denote it $\textbf{D}$.

\medskip

One can show that the following properties hold:

\medskip

\textbf{1.} For an arbitrary morphism $f\colon X \to Y$ of algebraic varieties we have
$$
\textbf{D}_Y \circ Rf_*\circ \textbf{D}_X \simeq  Rf_!
$$
and
\begin{align}\label{Eq.: Rf_*=DRf_!D}
Rf_*\simeq  \textbf{D}_Y \circ Rf_! \circ \textbf{D}_X.
\end{align}
One can be obtained from the other using $\textbf{D}\circ \textbf{D} \simeq \text{Id}$.

\medskip

\textbf{2.} If $g\colon Y \to \mathbb{A}^1$ is a morphism, then there exists a non-functorial isomorphism
\begin{align}\label{Eq.: D and Phi commute up to shift}
\textbf{D}\circ \Phi_g(\mathcal{F}^{\bullet}) \simeq \left( \Phi_g \circ \textbf{D} (\mathcal{F}^{\bullet})\right) [-2].
\end{align}

\textbf{3.} Duality commutes with restriction to an open subvariety.

\medskip
\textbf{4.} If $Y$ is smooth, then 
\begin{align}\label{Eq.: Dualizing complex for smooth variety}
\omega_Y\simeq \mathbb{C}_Y[2\dim Y].
\end{align}
Hence, $\textbf{D}_Y(\mathbb{C}_Y)\simeq\mathbb{C}_Y[2\dim Y]$. Here $\dim Y$ is the complex dimension.

\medskip

\textbf{5.} The functor $\textbf{D}$ preserves the support.

\subsection{Setup I.} 

Let $X$ be smooth complex algebraic variety, $U$ an open smooth subvariety, and $Z$ the complement to $U$. Let $\overline{f}\colon X \to \mathbb{A}^1$ be a proper morphism. Then we have a commutative diagram
\begin{align}\label{Eq.: MD}
\xymatrix{
U \ar[r]^j \ar[rd]_{f_U} &   X   \ar[d]^<<<<{\overline{f}}  & \ar[dl]^{f_Z} \ar[l]_i  Z         \\
    &        \mathbb{A}^1
}
\end{align}

(i) By \eqref{Eq.: Rf_*=DRf_!D} we have the isomorphism
$$
Rj_*(\mathbb{C}_U)\simeq  \textbf{D}_X \circ Rj_! \circ \textbf{D}_U(\mathbb{C}_U).
$$
Applying $\Phi_{\overline{f}}$ to both sides and using \eqref{Eq.: Dualizing complex for smooth variety} we get
$$
\Phi_{\overline{f}}\circ Rj_*(\mathbb{C}_U)\simeq  \Phi_{\overline{f}}\circ \textbf{D}_X \circ Rj_! (\mathbb{C}_U[2\dim U]).
$$

Using \eqref{Eq.: D and Phi commute up to shift} we can rewrite it as
$$
\Phi_{\overline{f}}\circ Rj_*(\mathbb{C}_U)\simeq   \left(\textbf{D}_X \circ \Phi_{\overline{f}}\circ Rj_! (\mathbb{C}_U[2\dim U])\right)[-2].
$$
If the support of the complex $\Phi_{\overline{f}}\circ Rj_! (\mathbb{C}_U)$ is  contained in $U$, then the same is true for $\Phi_{\overline{f}}\circ Rj_*(\mathbb{C}_U)$ (cf. Property 5 from Section \ref{SubSubSec.: Duality For Appendix}). Therefore, it does not matter which one to consider for cohomological tameness, since we are interested only in the support.

\medskip

(ii) There exists a standard short exact sequence (see \cite{Ha2}, Chapter II, exercise 1.19.c)
\begin{align}\notag
0 \to Rj_! \mathbb{C}_U \to \mathbb{C}_X \to i_* \mathbb{C}_Z \to 0.
\end{align}
Applying $\Phi_{\overline{f}}$ to it we get a distinguished triangle
\begin{align}\notag
\Phi_{\overline{f}} \circ Rj_! (\mathbb{C}_U) \to \Phi_{\overline{f}} (\mathbb{C}_X) \to \Phi_{\overline{f}} \circ i_* (\mathbb{C}_Z) \to \Phi_{\overline{f}} \circ Rj_! (\mathbb{C}_U)[1].
\end{align}
Since $i$ is \textit{proper}, we can rewrite it as
\begin{align}\notag
\Phi_{\overline{f}} \circ Rj_! (\mathbb{C}_U) \to \Phi_{\overline{f}} (\mathbb{C}_X) \to \hat{i}_* \circ \Phi_{{\overline{f}}\circ i} (\mathbb{C}_Z) \to \Phi_{\overline{f}} \circ Rj_! (\mathbb{C}_U)[1],
\end{align}
where $\hat{i}\colon Z_0 \to X_0$ is the natural closed embedding.

\smallskip

In the future applications we will need to prove that $\Phi_{\overline{f}} \circ Rj_! (\mathbb{C}_U)$ is supported on $U_0$. This would follow if we prove that  $\Phi_{\overline{f}} (\mathbb{C}_X)$ is supported on $U_0$ and $\Phi_{\overline{f}\circ i} (\mathbb{C}_Z)$ has empty support (since $\hat{i}_*$ does not change support).

Complexes $\Phi_{\overline{f}} (\mathbb{C}_X)$ and $\Phi_{{\overline{f}}\circ i} (\mathbb{C}_Z)$ compute vanishing cycles of $\overline{f}$ and $f_Z$ with values in $\mathbb{C}_X$ and $\mathbb{C}_Z$ respectively. Therefore, their support can be computed geometrically.

\subsection{Setup II.}

It is clear that proving cohomological tameness of \eqref{Eq.: F for Q_3} is equivalent to proving cohomological tameness of 
\begin{align}\notag
g(x,y,z)=y(z+1)+\frac{qx^2}{(xy-1)z}.
\end{align}
Here $x,y,z$ are coordinates on $\bA^3=\Spec (\bC[x,y,z])$ and $g$ is defined on $U \subset \bA^3$ given by the equation $(xy-1)z \neq 0$. By a simple computation of partial derivatives one can see that all critical points of $g$ lie in the subset $\{xy-z-1=0\}$.

Consider the map $U \to U\times \bA^1_t$ given by the graph of $g$ and denote its image $\Gamma_g \subset U\times \bA^1_t$. Further, consider the natural embedding of $U\times \bA^1_t$ into  $\bA^3 \times \bA^1_t \subset (\bP^1)^3 \times \bA^1_t$. Let $\overline{\Gamma}_g$ be the closure of $\Gamma_g$ in $(\bP^1)^3 \times \bA^1_t$. Thus, we have a commutative diagram completely analogous to (\ref{Eq.: MD})
\begin{align}
\xymatrix{
\Gamma_g \ar[r]^j \ar[rd]_{g_U} &  \overline{\Gamma}_g   \ar[d]^<<<<{\overline{g}}  & \ar[dl]^{g_Z} \ar[l]_i  \overline{\Gamma}_g \setminus \Gamma_g         \\
    &        \mathbb{A}^1_t
}
\end{align}
where $\overline{g}$ is induced by projection $(\bP^1)^3 \times \bA^1_t\to \bA^1_t$, and we will identify $g_U$ with $g$.

It is easy to see that $\Gamma_g \subset \bA^3 \times \bA^1_t$ is defined by the equation
\begin{align}\notag
y(z+1)+\frac{qx^2}{(xy-1)z}=t,
\end{align}
and $\overline{\Gamma}_g \subset (\mathbb{P}^1)^3 \times  \bA^1_t$ is given by the homogeneous equation
\begin{align} \label{Eq.: HomogEquationForG_f}
x_0z_1\left(x_1 y_1-x_0y_0\right)\left[y_1\left(z_1+z_0\right)-ty_0z_0\right]+qz_0^2x_1^2y_0^2=0.
\end{align}

\subsection{Open cover.} 

Consider an open cover of $(\mathbb{P}^1)^3 \times \bA^1_t$ by 8 open subsets
$$
V_{i,j,k}=\{x_iy_jz_k\neq 0\},
$$
each of which is just $\bA^3\times \bA^1_t$. As standard local coordinates on these open subsets we will be always
using fractions $\frac{x_{i+1}}{x_i}$, $\frac{y_{j+1}}{y_j}$, $\frac{z_{k+1}}{z_k}$ and $t$ (in the subscripts here we mean mod 2 sum: $1+1=0$).

On each $V_{i,j,k}$ one can write down the equation of $\overline{\Gamma}_g \cap V_{i,j,k}$ in terms of local coordinates and we will be referring to this equation simply as "the equation of $\overline{\Gamma}_g$ in the chart $V_{i,j,k}$".

\subsection{How to compute $\Phi_{\overline{g}} (\mathbb{C}_{\overline{\Gamma}_g})$?} 

Decompose $\overline{\Gamma}_g$ into the disjoint union
$$
\overline{\Gamma}_g=W \sqcup S,
$$
where $S$ is the singular locus of $\overline{\Gamma}_g$ and $W$ is the smooth locus; $S$ is a closed subvariety of codimension at least 1, $W$ is open in $\overline{\Gamma}_g$ and smooth.

According to Section \ref{SubSubSec.: RestrictionToOpen}, computing vanishing cycles commutes with restriction to an open subset, and therefore, to investigate the support at points of $W$, we can restrict to it from the beginning. Since $W$ is smooth the support of $\Phi_{\overline{g}_{|W}} (\mathbb{C}_W)$ can be non-zero only in singular points of fibers of $\overline{g}_{|W}$ (by the implicit function theorem). Thus, on $W$ we just need to worry about critical points of $\ol{g}_{|W}$.

\smallskip

There exists a decomposition into disjoint union
$$
(\mathbb{P}^1)^3\times \bA^1_t= V_{0,0,0} \sqcup \{x_0y_0z_0=0 \}.
$$

\subsubsection{Lemma.}

{\it Let $Q \in \overline{\Gamma}_g$ be a smooth point, i.e. $Q \in W$. If $Q \in \{x_0y_0z_0=0 \}$, then $Q$ is not a critical point of $\ol{g}_{|W}$.
}

\smallskip

\proof Consider a chart $V_{i,j,k}$ and temporarily denote the local coordinates just by $x,y,z$. Let $P(x,y,z,t)=0$ be equation (\ref{Eq.: HomogEquationForG_f}) written in these coordinates. 

\smallskip

In this notation the intersection $V_{i,j,k} \cap S$ is defined by the system
\begin{align}\label{Eq.: Singular points for bar(gamma g)}
&x_0z_1(x_1 y_1-x_0y_0)y_0z_0=0 \\ \notag
&P_x(x,y,z,t)=0 \\ \notag
&P_y(x,y,z,t)=0 \\ \notag 
&P_z(x,y,z,t)=0 \\ \notag
&P(x,y,z,t)=0,
\end{align}
where we have written the first equation, which comes from the derivative with respect to $t$, in the original homogeneous coordinates.

\smallskip

On the other hand, the intersection of the fiber of $\ol{g}$ over the point $d \in \bA^1_t$ with $V_{i,j,k}$ is given by the equation $P(x,y,z,d)=0$. Therefore, on $V_{i,j,k}$ the singular locus of this fiber is given by the system
\begin{align}\label{Eq.: Singular points for a fiber of bar(g)}
&P_x(x,y,z,d)=0 \\  \notag
&P_y(x,y,z,d)=0 \\  \notag
&P_z(x,y,z,d)=0 \\  \notag
&P(x,y,z,d)=0.
\end{align}

From \eqref{Eq.: HomogEquationForG_f} it is easy to see that if a point $Q\in \{x_0y_0z_0=0 \}$ with coordinates $(a,b,c,d)$ satisfies \eqref{Eq.: Singular points for a fiber of bar(g)}, then it satisfies it for arbitrary $d$, i.e. we get singular points simultaneously in all fibers. Thus, on the locus $\{x_0y_0z_0=0 \}$ systems \eqref{Eq.: Singular points for a fiber of bar(g)} and \eqref{Eq.: Singular points for bar(gamma g)} coincide.

\smallskip

Therefore, if $Q \in \{x_0y_0z_0=0 \}$ is a smooth point of $\ol{\Gamma}_g$, then it is also a smooth point in the respective fiber of $\ol{g}$ (and $\ol{g}_{|W}$). $\blacksquare$

\subsubsection{Strategy.}

In the rest of this section we will treat points in  $V_{0,0,0}$ and in $\{x_0y_0z_0=0 \}$ separately. In the latter case, by the above lemma, it is enough to look at vanishing cycles at singular points of $\ol{\Gamma}_g$.

\subsection{Chart $V_{0,0,0}$.} 

Recall that $\Gamma_g \subset V_{0,0,0}$ and in this chart we need to prove that the support is contained in $\Gamma_g$. In this chart $\overline{\Gamma}_g$ is given by
$$
z(xy-1)[y(1+z)-t]+qx^2 = 0
$$
and $\Gamma_g$ is defined by the intersection with $z(xy-1) \neq 0$. Consider functions
\begin{align}  \label{Eq.: NewCoord000}
&x_1=x\\ \notag
&y_1=y(1+z)-t\\ \notag
&z_1=z(xy-1)\\ \notag
&t_1=t.
\end{align}
Computing the Jacobian we get
$$
J=xy-z-1,
$$
and hence on the complement to the closed subset $\{xy-z-1=0\}$ formulas (\ref{Eq.: NewCoord000}) define a new coordinate system. Notice that $\left(\overline{\Gamma}_g \setminus \Gamma_g\right) \cap V_{0,0,0}=\{x=0,\, z=0\}$ lies in this open set.

On this open set the equation for $\overline{\Gamma}_g$ can be rewritten as
\begin{align}\label{Eq.: V000 new equation}
y_1z_1+qx^2=0,
\end{align}
and hence the restriction of $\overline{g}$ to this open set is a projection and therefore has no vanishing cycles. Thus, we conclude that in this chart all vanishing cycles of $\overline{g}$ with coefficients in $\mathbb{C}_{\overline{\Gamma}_g }$ live in $\Gamma_g$.

Using \eqref{Eq.: V000 new equation} it is easy to see that 
$$
V_{0,0,0}\cap S=\{x=0, \, z=0, \, y=t \}.
$$

\subsection{The set  $\{x_0y_0z_0=0 \}$.} 

To check vanishing of stalks of $\Phi_{\overline{g}} (\mathbb{C}_{\ol{\Gamma}_g})$ at points of $S$ we will be restricting to different charts and work in local coordinates. The whole chart $V_{0,0,0}$ has already been considered. Hence, there are seven charts left to be checked.

\subsubsection{Charts $V_{1,0,0}, V_{0,1,0}, V_{0,0,1}$.}

In the chart $V_{0,1,0}$ the equation of $\overline{\Gamma}_g$ takes the form
\begin{align} \label{Eq.: olGfV010}
z(x-y')[(z+1)-y't]+qx^2y'^2=0.
\end{align}
Finding singular points of $\overline{\Gamma}_g$ with $y'=0$ (we need to check only them) we get that there are two singular lines given by
\begin{align} \label{Eq.: SingV010}
&x=0, \, y'=0, \, z=0  \\ \notag
&x=0, \, y'=0, \, z=-1.
\end{align}
By computing the Jacobian it is easy to see that functions
\begin{align}\label{Eq.: NewCoord010}
x_1=x, \quad y_1=y', \quad z_1=(1+z-ty')z, \quad t_1=t.
\end{align}
form a coordinate system in a neighborhood of (\ref{Eq.: SingV010}). In terms of these coordinates (\ref{Eq.: olGfV010}) takes the form
\begin{align}\notag
z_1(x_1-y_1)+qx_1^2y_1^2=0,
\end{align}
and is independent of $t$. Hence, there are no vanishing cycles.

\medskip

It is easy to check that the sets $S\cap V_{1,0,0} \cap \{x_0=0 \}$ and $S\cap V_{0,0,1} \cap \{z_0=0 \}$ are empty.

\subsubsection{Charts $V_{1,1,0}, V_{0,1,1}$.} 

In the previous section we have considered the charts $V_{1,0,0}, V_{0,1,0}, V_{0,0,1}$. Therefore, in the chart $V_{1,1,0}$ we only need to look at points in $S\cap V_{1,1,0} \cap \{x_0=0 \}\cap \{y_0=0 \}$. Keeping in mind this remark one can treat the chart $V_{1,1,0}$ completely similarly to $V_{0,1,0}$ just replacing $x$ by $x'$ everywhere; the conclusion is also similar. 

An analogous remark applies to $V_{1,0,1}$ and $V_{0,1,1}$ as well. Looking for points in $S\cap V_{0,1,1} \cap \{y_0=0 \}\cap \{z_0=0 \}$ it is easy to see that this set is empty. The chart $V_{1,0,1}$ will be considered in the next section.

\subsubsection{Chart $V_{1,0,1}$.} 

This is the only case where the argument is a bit more involved. In this chart the equation of $\overline{\Gamma}_g$ takes form
\begin{align}\notag
x'(y-x')[y(1+z')-tz']+qz'^2=0,
\end{align}
and the singular locus with $x'=0$ and $z'=0$ is the line
$$
x'=0, \, y=0, \, z'=0.
$$
It is not clear if one can find such a simple argument via coordinate changes as above. Therefore, we proceed differently.

For each value of $t$ the fiber over it has an isolated singularity at the origin, and we need to prove that there are no vanishing cycles to this point. Note that the family is of the form 
$$
\alpha(x', y, z') + t \beta(x', y, z')=0,
$$
where
\begin{align*}
&\alpha(x', y, z')=x'(y-x')y(1+z')+qz'^2 \\
&\beta(x', y, z')=-x'(y-x')z'.
\end{align*}
According to Corollary 2.1 of \cite{Pa}, if this family of isolated hypersurface singularities is $\mu$-constant (i.e. Milnor numbers coincide for all fibers), then it is topologically locally trivial over the base. In particular, this implies the absence of vanishing cycles.

Let us compute the Milnor numbers. Consider the polynomial $P(x', y, z')=x'(y-x')[y(1+z')-tz']+qz'^2$. Its partial derivatives are
\begin{align} \notag
&P_{x'}=(y-2x')[y(1+z')-tz']\\ \notag
&P_{y}=x'[y(1+z')-tz']+x'(y-x')(1+z')=x'[(2y-x')(1+z')-tz'] \\ \label{Eq.: Derivative wrt z'}
&P_{z'}=x'(y-x')[y-t]+2qz'.
\end{align}

By definition the Milnor number is 
\begin{align*}
\mu=\dim_{\bC} \left( \bC\{x', y, z'\}/(P_{x'}, P_y, P_{z'}) \right),
\end{align*}
and we need to show that it does not depend on $t$.

From \eqref{Eq.: Derivative wrt z'} we see that $2qz'= - x'(y-x')[y-t]$ in the quotient $\bC\{x', y, z'\}/(P_{x'}, P_y, P_{z'})$. Hence, the latter can be rewritten as
\begin{align}\label{Eq.: Milnor Ring I}
\bC\{x', y\}/(\wt{P}_{x'}, \wt{P}_y), 
\end{align}
where
\begin{align*}
&\wt{P}_{x'}= (y-2x')[y(1-\frac{x'(y-x')[y-t]}{2q})+t\frac{x'(y-x')[y-t]}{2q}]\\
&\wt{P}_y=x'[(2y-x')(1-\frac{x'(y-x')[y-t]}{2q})+t\frac{x'(y-x')[y-t]}{2q}].
\end{align*}
Note that we can write
\begin{align*}
&\wt{P}_{x'}=f_1 f_2, \quad \wt{P}_y=f_3f_4,
\end{align*}
where
\begin{align*}
&f_1= (y-2x'), \quad \quad f_2=[y(1-\frac{x'(y-x')[y-t]}{2q})+t\frac{x'(y-x')[y-t]}{2q}]\\
&f_3=x', \quad \quad f_4=[(2y-x')(1-\frac{x'(y-x')[y-t]}{2q})+t\frac{x'(y-x')[y-t]}{2q}].
\end{align*}
It easy to check that any pair out of $f_1,f_2, f_3 , f_4$ forms a coordinate system around the origin in  the $x', y$-plane. Let $u=f_1, v=f_2$ be such a coordinate system. Then \eqref{Eq.: Milnor Ring I} can be rewritten as 
\begin{align}\label{Eq.: Milnor Ring II}
\bC\{u, v\}/(uv, \wt{P}_y), 
\end{align}
where $\wt{P}_y$ is written as a power series in $u$ and $v$; it starts from terms quadratic in $u$ and $v$, since $\wt{P}_y=f_3f_4$.  Moreover, as a vector space \eqref{Eq.: Milnor Ring II} can be further rewritten as
\begin{align*}
\bC\{u\}/(u^2) \oplus \bC\{v\}/(v^2).
\end{align*}
Thus, 
\begin{align*}
\mu=4,
\end{align*}
and is independent of $t$.

\subsubsection{Chart $V_{1,1,1}$.}

It is easy to check that the set $S\cap V_{1,1,1} \cap \{x_0=0 \} \cap \{y_0=0 \} \cap \{z_0=0 \}$ is empty.

\subsection{Computation of $\Phi_{\overline{f}\circ i} (\mathbb{C}_Z)$.} 

Recall that $Z=\overline{\Gamma}_g\setminus \Gamma_g$ and $\overline{\Gamma}_g \subset (\mathbb{P}^1)^3 \times \mathbb{C}$ is given by the homogeneous equation
\begin{align}\notag
x_0z_1\left(x_1 y_1-x_0y_0\right)\left[y_1\left(z_1+z_0\right)-ty_0z_0\right]+qz_0^2x_1^2y_0^2=0.
\end{align}
The subvariety $\Gamma_g$ is defined by additionally putting
\begin{align}\notag
&x_0y_0z_0\neq 0 \\ \notag
&z_1(x_1y_1-x_0y_0)\neq 0.
\end{align}
Thus, $\overline{\Gamma}_g\setminus \Gamma_g$ is defined by the system
\begin{align}\notag
&x_0z_1\left(x_1 y_1-x_0y_0\right)\left[y_1\left(z_1+z_0\right)-ty_0z_0\right]+qz_0^2x_1^2y_0^2=0 \\ \notag
&x_0y_0z_0 z_1(x_1y_1-x_0y_0)= 0.
\end{align}
This system is equivalent to
\begin{align}\notag
&x_0z_1y_1(x_1 y_1-x_0y_0)(z_1+z_0)+qz_0^2x_1^2y_0^2=0 \\ \notag
&x_0y_0z_0 z_1(x_1y_1-x_0y_0)= 0,
\end{align}
which is independent of $t$. Therefore $Z=\overline{\Gamma}_g\setminus \Gamma_g$ is a product and $\Phi_{\overline{g}\circ i} (\mathbb{C}_Z)$ has empty support.

{\small
\textsc{Vassily Gorbounov, Institute of Mathematics, University of Aberdeen, Aberdeen, AB24 3UE, UK}

\textit{E-mail address}: \texttt{v.gorbunov@abdn.ac.uk}

\smallskip

\textsc{Maxim Smirnov, Mathematics Section, The Abdus Salam International Centre for Theoretical Physics, Strada Costiera 11, 34151 Trieste, Italy}

\textit{E-mail address}: \texttt{msmirnov@ictp.it}

}

\end{document}